\documentclass[10pt]{article}
\usepackage{amssymb}
\usepackage{enumerate}
\usepackage{graphicx}
\usepackage{amsthm}
\usepackage{color}

\def\R{{\rm I\! R}}

\def\sss{strictly star-shaped }
\def\fx{{\partial \phi \over \partial x}}
\def\fy{{\partial \phi \over \partial y}}
\def\px{{\partial P \over \partial x}}
\def\py{{\partial P \over \partial y}}
\def\qx{{\partial Q \over \partial x}}
\def\qy{{\partial Q \over \partial y}}
\def\hkx{{\partial H_{2k} \over \partial x}}
\def\hky{{\partial H_{2k} \over \partial y}}

\newtheorem{theorem}{Theorem}

\newtheorem{lemma}{Lemma}

\title{Existence and uniqueness of limit cycles in a class of second order ODE's }

\author{M. Sabatini  \footnote{Dip. di Matematica, Univ. di Trento, I-38050 Povo, (TN) - Italy.
Email: marco.sabatini@unitn.it,
Phone: ++39(0461)881670, Fax: ++39(0461)881624 }}
\date{March 3, 2010}
\begin{document}
\maketitle
\begin{abstract}We prove a uniqueness result for limit cycles of a class of second order ODE's. As a special case, we prove limit cycle's uniqueness for an ODE studied in \cite{ETBA}.
\end{abstract}

\section{Introduction}
 
Let us consider a first order differential system in the real plane,
\begin{equation}\label{sysPQ}
\dot x = P(x,y), \qquad \dot y = Q(x,y).
\end{equation} 
The study of the dynamics of (\ref{sysPQ}) strongly depends on the existence and stability properties of special solutions such as equilibrium points and non-constant periodic solutions. In particular, if an attracting non-constant periodic solution exists, then it dominates the dynamics of (\ref{sysPQ})  in an open, connected subset of the plane, its region of attraction. In some cases such a region of attraction can even extend to cover the whole plane, with the unique exception of an equilibrium point. Uniqueness theorems for non-constant periodic solutions, i. e. limit cycles, have been  extensively studied,  see \cite{CRV} and \cite{XZ} for recent results and extensive bibliographies. Most of the results known are concerned with the classical Li\'enard system, 
\begin{equation}\label{syslie} 
\dot x = y - F(x)   , \qquad \dot y = - g(x). 
\end{equation} 
and its generalizations, such as
\begin{equation}\label{sysGG}
\dot x = \beta(x)\big[ \varphi(y) - F(x)  \big], \qquad \dot y = -\alpha(y)g(x). 
\end{equation}
Such a class of systems also contain Lotka-Volterra systems and systems equivalent to Rayleigh equation
\begin{equation}\label{equaray} 
\ddot x + f(\dot x) + g(x) = 0,
\end{equation} 
as special cases. A very recent result \cite{CRV} is concerned with systems equivalent to
\begin{equation}\label{equaCRV} 
\ddot x +  \sum_{k=0}^{N}f_{2k+1}(x){\dot x}^{2k+1} + x = 0,
\end{equation} 
with $f_{2k+1}(x)\geq 0$, increasing  for $x > 0$, decreasing for $x < 0$, $k=0, \dots, N$.
On the other hand, there exist classes of second order ODE's which are not covered by the above cases. This is the case of a model developped in \cite{ETBA}, which led to the equation
\begin{equation}\label{ETBA}
\ddot x + \epsilon \dot x (x^2 + x \dot x + {\dot x}^2 -1) + x = 0, \qquad \epsilon >0.
\end{equation}
In this paper we prove a uniqueness result for systems equivalent to
\begin{equation}\label{equaphi}
\ddot x + \dot x \phi(x,\dot x) + x = 0,
\end{equation}
under the assumtpion that $\phi(x,y)$ be a  function with star-shaped level sets. As a consequence, we are able to prove existence and uniqueness of the limit cycle for the equation (\ref{ETBA}).

\section{Risultati preliminari}

Let $\Omega \subset \R^2$ be a star-shaped set. We say that a function $\phi \in C^1( \Omega, \R) $ is  {\it star-shaped} if $(x,y) \cdot \nabla \phi = x \fx + y \fy$ does not change sign. We say that $\phi$ is 
{\it strictly star-shaped} if  $(x,y) \cdot \nabla \phi \neq 0$.   We call {\it ray} a half-line having origin at the point $(0,0)$.  

Let us consider a system equivalent to the equation (\ref{equaphi})
\begin{equation}\label{sisphi}
\dot x = y \qquad \dot y = -x - y \phi(x,y).
\end{equation}
We denote by $\gamma(t,x^*,y^*)$ the unique  solution to the system (\ref{sisphi}) such that $\gamma(0,x^*,y^*) = (x^*,y^*)$. We first consider a sufficient condition for limit cycles' uniqueness.

\begin{theorem}\label{teorema} Let $\phi: \R^2 \rightarrow \R^2 $ be a \sss function. Then (\ref{sisphi}) has at most one limit cycle.
\end{theorem}
{\it Proof.}
Let us assume that, for $(x,y) \neq  (0,0)$, 
$$
 x \fx + y  \fy > 0.
$$
The proof can be performed analogously for the opposite inequality.
 
Applying Corollary 6 in \cite{GS} requires to  compute the expression
$$
\nu = P\left( x\qx + y\qy \right) - Q \left( x\px + y\py\right) ,
$$
where $P$ and $Q$ are the components of the considered vector field. For system (\ref{sisphi}), one has
$$
\nu =
y \left(-x - xy \fx - y\phi - y ^2 \fy \right) - \left(  -x - y \phi(x,y) \right) y  =
$$
$$
-y^2 \left( x \fx + y  \fy \right) \leq 0.
$$
The function $\nu$ vanishes only for $y=0$. Let us assume, by absurd, that two distinct limit cycles exist, $\gamma_1$ and $\gamma_2$. Since the system  (\ref{sisphi}) has only one critical point, the two cycles have to be concentric. Let us assume that  $\gamma_2$ encloses $\gamma_1$. For both cycles one has:
$$
\int_{0}^{T_i} \nu(\gamma_i(t) )dt < 0, \qquad i=1,2,
$$
where $T_i$ is the period of $\gamma_i$, $i=1,2$. Hence both cycles, by theorem 1 in \cite{GS}, are attractive. Let $A_1$ be the region of attraction of $\gamma_1$.  $A_1$ is bounded, because it is enclosed by $\gamma_2$, which is not attracted to $\gamma_1$. The external component of   $A_1$'s boundary is itself a cycle $\gamma_3$, because (\ref{sisphi}) has just one critical point at the origin. Again,
$$
\int_0^{T_3} \nu(\gamma_3(t) )dt < 0,
$$
hence $\gamma_3$ is attractive, too. This contradicts the fact that the solutions of (\ref{sisphi}) starting from its inner side are attracted to $\gamma_1$. Hence the system  (\ref{sisphi}) can have at most a single limit cycle.
\hfill $\clubsuit$

In particular, the equation (\ref{ETBA}) considered in \cite{ETBA} has at most one limit cycle. In fact, in this case one has $\phi(x,y) = \epsilon (x^2 + xy + y^2 -1)$, so that one has 
$$
\nu = x \fx + y  \fy = 2 \epsilon y^2 (x^2 + xy + y^2) > 0  \quad {\rm for } \quad  (x,y) \neq (0,0).
$$
It should be noted that even if the proof is essentially based on a stability argument, the divergence cannot be used in order to replace the function $\nu$. In fact, the divergence of system (\ref{sisphi}) is
$$
{\rm div } \big(y , -x - y \phi(x,y) \big)= - \phi -  y  \fy, 
$$
which does not have constant sign, under our assumptions. Moreover, the divergence cannot have constant sign in presence of a repelling critical point and an attracting cycle.

Now we care about the existence of limit cycles. We say that $\gamma(t)$ is {\it positively bounded} if the semi-orbit $\gamma^+ = \{\gamma(t), \quad t \geq 0\}$ is contained in a bounded set. Let us denote by $D_r$ the disk $ \{ (x,y) : dist((x,y),O) \leq r \} $, and by $B_r$ its boundary $\{ (x,y) : dist((x,y),O) = r \} $. 
In the following, we use the function $V(x,y) =  \frac {x^2}2 + \frac {y^2}2$ as a Liapunov function.
Its derivative along the solutions of (\ref{sisphi}) is
$$
\dot V(x,y) = - y^2 \phi(x,y).
$$

\begin{lemma}\label{lemma} Let $U$ be a bounded set, with $\sigma := \sup \{  dist((x,y),O), (x,y) \in U \}$. 
If $\phi(x,y) \geq 0$ out of $U$, and $\phi(x,y)$ does not vanish identically on any $B_r$, for $r > \sigma$, then every $\gamma(t)$ definitely enters the disk $D_\sigma$ and does not leave it.
\end{lemma}
{\it Proof.}
The level curves of $V(x,y)$ are circumferences. For every $r \geq \sigma$, the disk $D_r $ contains $U$.  Since $\dot V(x,y) = - y^2 \phi(x,y) \leq 0$ on its boundary, such a disk is positively invariant. 
Let $\gamma$ be an orbit with a point $\gamma(t^*)$ such that $d^* = dist(\gamma(t^*),O) > \sigma$. Then $\gamma$ does not leave the disk $D_{d^*}$, hence it is positively bounded.
Moreover $\gamma(t)$ cannot be definitely contained in $B_r$, for any $r > \sigma$, since $\dot V(x,y)$ does not vanish identically on any $B_r$, for $r > \sigma$. Now, assume by absurd that $\gamma(t)$ does not intersect  $B_\sigma$. Then its positive limit set is a cycle $\overline{\gamma}(t)$, having no points in $D_\sigma$. The cycle $\overline{\gamma}(t)$ cannot cross outwards any $B_r$, hence it has to be contained in $B_r$, for some $r > \sigma$, contradicting the fact that $\dot V(x,y)$ does not vanish identically on any $B_r$, for $r > \sigma$. Hence there exists  $t^+ > t^*$ such that $\gamma(t^+) \in D_\sigma$. Then, for every  $t > t^+$, one has $\gamma(t) \in D_\sigma$, because $\dot V(x,y) \leq 0$ on $B_\sigma$.
\hfill$\clubsuit$

Collecting the results of the above statements, we may state a theorem of existence and uniqueness for limit cycles of a class of second order equations. We say that an equilibrium point $O$ is {\it negatively asymptotically stable} if it is asymptotically stable for the system obtained by reversing the time direction.

\begin{theorem} If the hypotheses of theorem \ref{teorema} and lemma \ref{lemma} hold, and $\phi(0,0) < 0$, then the system (\ref{sisphi}) has exactly one limit cycle, which attracts every non-constant solution.
\end{theorem}
{\it Proof.} By the above lemma, all the solutions are definitely contained in $D_\sigma$. The condition $\phi(0,0) < 0$  implies by continuity $\phi(x,y) < 0$ in a neighbourhood $N_O$ of the origin. This gives  the negative asymptotic stability of the origin by Lasalle's invariance principle, since $\dot V(x,y) \geq 0$ in $N_O$, and the set $\{\dot V(x,y) = 0\} \cap N_O = \{y = 0\}   \cap N_O$ does not contain any positive semi-orbit. The system has just one critical point at the origin, hence by Poincar\'e-Bendixson theorem there exist a limit cycle. By  theorem \ref{teorema}, such a limit cycle is unique.
\hfill$\clubsuit$

This proves that every non-constant solution to the equation (\ref{ETBA}) studied in  \cite{ETBA} is attracted to the unique limit cycle. 

We can produce more complex systems with such a property. Let us set
$$
\phi(x,y) = -M +\sum_{k=1}^n H_{2k}(x,y),
$$
with $ H_{2k}(x,y)$ is a homogeneous function of degree $2k$, positive except at the origin, $M$ is a positive constant. Then, by Euler's identity,  one has
$$
\nu = \sum_{k=1}^n \left( x\hkx + y\hky \right) = \sum_{k=0}^n 2kH_{2k}(x,y) > 0  \quad {\rm for } \quad  (x,y) \neq (0,0).
$$
If $\phi(x,y)$ does not vanish identically on any $B_r$, for instance if $H_{2k}(x,y) = (x^2 +xy + y^2)^k$, then the corresponding system (\ref{sisphi}) has a unique limit cycle.
In general, it is not necessary to assume the positiveness of all of the homogeneous functions $H_{2k}(x,y)$, as the following example shows. Let us set $Q(x,y) = x^2+xy+y^2$. Then take
$$
\phi(x,y) =  -1 + Q  - Q^2 + Q^3 .
$$
One has
$$
\nu = x \fx + y  \fy = 2Q - 4Q^2 +6Q^3 = Q(2  - 4Q  +6Q^2)  .
$$
The discriminant of the quadratic polynomial $2  - 4Q  +6Q^2$ is $\Delta =  -32 < 0$ hence $\nu >0$ everywhere but at the origin. Moreover, $\phi(x,y)$ does not vanish identically on any circumference, hence  the corresponding system (\ref{sisphi}) has a unique limit cycle.

\end{document}